\documentclass{article}
\usepackage{latexsym}
\usepackage{amssymb}
\begin{document}
\def\P{\Bbb P}
\def\Q{\Bbb Q}
\title{Morasses and Finite Support Iterations}
\author{Bernhard Irrgang}
\date{}
\maketitle
\begin{abstract}
We introduce a method of constructing a forcing along a simplified $(\kappa,1)$-morass such that the forcing satisfies the $\kappa$-chain condition. Alternatively, this may be seen as a method to thin out a larger forcing to get a chain condition. As an application, we construct a ccc forcing that adds an $\omega_2$-Suslin tree. Related methods are Shelah's historic forcing and Todorcevic's $\rho$-functions.
\end{abstract}
\section{Introduction}
There are a number of consistency questions from two-cardinal combinatorics that were answered by Shelah's method of historic forcing or with the help of Todorcevic's $\rho$-functions: Can there exist a superatomic Boolean algebra with width $\omega$ and height $\omega_2$ (Baumgartner and Shelah \cite{BaumgartnerShelah}, Martinez  \cite{Martinez})?
Is it possible, that there is a function $f:\omega_2 \times \omega_2 \rightarrow \omega$, such that $f$ is non-constant on any rectangle with infinite sides (Todorcevic  \cite{Todorcevic,Todorcevic2})?
Can one prove in $ZFC$, that every initially $\omega_1$-compact $T_3$-space with countable tightness is already compact (Rabus  \cite{Rabus}, Juhasz and Soukup \cite{JuhaszSoukup})?
Is there consistently a forcing that satisfies ccc and adds a Kurepa tree (Jensen \cite{Jensen2}, Velickovic \cite{Velickovic})? There are many more examples, but we cannot give a comprehensive overview here. In all cases there is a natural forcing with finite conditions that would solve the problem if it preserved cardinals. Since the conditions are finite, the suitable property of the forcing to guarantee cardinal preservation is the countable chain condition (ccc). Therefore one thins out the natural forcing in such a way that the remaining forcing satisfies ccc.
\smallskip\\
In the following, we will present the simplest case of a morass approach to such questions, i.e. to construct a ccc-forcing of size $\omega_2$. The basic idea is simple: We try to generalize iterated forcing with finite support. Classical iterated forcing with finite support as introduced by Solovay and Tennenbaum \cite{SolovayTennenbaum} works with continuous, commutative systems of complete embeddings of Boolean algebras or partial orders which are indexed along a well-order. The following holds: If every forcing of the system 
satisfies ccc, then also the direct limit does. So if e.g. all forcings of the system are countable, then its direct limit satisfies ccc. It will, however, have size $\leq \omega_1$ since it is a direct limit, while we want to construct a forcing of size $\omega_2$. To overcome this limitation, we will not consider a linear system indexed along a well-order but a two-dimensional system indexed along a simplified $(\omega _1,1)$-morass. Since we want to obtain complete embeddings, we have to thin out the natural forcings. The way to do this follows very naturally from our approach. As an example how the thinning out is done, we will construct a ccc forcing that adds an $\omega_2$-Suslin tree. The basic forcing we thin out is Tennenbaum's forcing for adding a Suslin tree with finite conditions.
\smallskip\\
Morasses were introduced by R. Jensen in the early 1970's to solve the cardinal transfer problem of model theory in $L$ (see e.g. Devlin \cite{Devlin}). For the proof of the gap-2 transfer theorem a gap-1 morass is used. For higher-gap transfer theorems Jensen has developed so-called higher-gap morasses \cite{Jensen1}.  In his Ph.D. thesis, the author generalized these to gaps of arbitrary size (see \cite{Irrgang3,Irrgang1,Irrgang2}). The theory of morasses is far developed and well examined. In particular, it is known how to construct morasses in $L$ \cite{Devlin,Friedman,Irrgang3,Irrgang2} and how to force them \cite{Stanley2,Stanley}. Moreover, D. Velleman has defined so-called simplified morasses, along which morass constructions can be carried out more easily \cite{Velleman1984,Velleman1987a,Velleman1987b}. Their existence is equivalent to the existence of ordinary morasses \cite{Donder, Morgan1998}. The fact that the theory of morasses is so far developed is an advantage of the morass approach compared to historic forcing or $\rho$-functions. It allows straightforward generalizations to higher cardinals while the conditions of the forcings can be kept finite.
\smallskip\\ 
While the general method presented here works for higher cardinals, we can in general not expect that the consistency statements can naively be extended by raising the cardinal parameters. For example, we force an $\omega_2$-Suslin tree along a gap-1 morass. An innocent generalization of the argument that the resulting tree has neither a branch nor an antichain of size $\omega_2$, would yield a tree on $\omega_3$ that has neither a branch nor an antichain of size $\omega_2$, which is of course impossible. The reason why this generalization does not work is that the gap-2 case yields a three-dimensional construction. Therefore, the finite conditions of our forcing have to fit together appropriately in three directions instead of two directions and that is impossible. So if and how a statement generalizes to higher gaps depends  heavily on the concrete conditions.
\smallskip\\
The exact relationship between our approach and the methods of historic forcing and $\rho$-functions is an open question. The crucial step in our proof that chain conditions are preserved is the definition of the support of a condition. It resembles the definition of the ``history'' $t^\ast(\alpha)$ of an ordinal $\alpha$ given by Baumgartner and Shelah \cite{BaumgartnerShelah}. However, there are various ways to set things up and the definition of a FS system given below is just one of them. As far as $\rho$-functions are concerned, it is possible to directly read off a $\rho$-function from a simplified gap-1 morass. This is a result of C. Morgan's in \cite{Morgan1996}. It is, however, unclear how this relates to an approach as below which generalizes finite support iterations.
\smallskip\\
If $\P$ is the limit of a finite support iteration indexed along $\alpha$, then we can understand a $\P$-generic extension as being obtained successively in $\alpha$-many steps. Moreover, there are names for the forcings used in every step. This raises the question if a similar analysis is possible for a forcing which is constructed with our method. It would justify to call them FS iterations along morasses instead of FS systems along a morass, which was the name the author used until the referee pointed out the shortcoming concerning successive extensions.
\smallskip\\
We should also mention that besides historic forcing and $\rho$-functions there is another, quite different method to prove consistencies in two-cardinal combinatorics. This is the method of forcing with models as side conditions or with side conditions in morasses. Models as side conditions were introduced by S. Todorcevic \cite{Todorcevic1985,Todorcevic1989}, which was further developed by P. Koszmider \cite{Koszmider2000} to side conditions in morasses. Unlike the other methods, it produces proper forcings which are usually not ccc. This is sometimes necessary. For example, Koszmider proved that if CH holds, then there is no ccc forcing that adds a sequence of $\omega_2$ many functions $f:\omega_1 \rightarrow \omega_1$ which is ordered by strict domination mod finite. However, he is able to produce a proper forcing which adds such a sequence \cite{Koszmider2000}. More on the method can be found in Morgan's paper \cite{Morgan2006}. In the context of our approach, this raises the question if it is possible to define something like a countable support iteration along a morass. 
\section{FS Iterations}
Let $\P$ and $\Q$ be partial orders. A map $\sigma :\P \rightarrow \Q$ is called a complete embedding if
\smallskip\\
(1) $\forall  p ,p^\prime \in \P \ (p^\prime \leq p \rightarrow \sigma (p^\prime ) \leq \sigma (p))$
\smallskip\\
(2) $\forall  p ,p^\prime \in \P \ ( p$ and $p^\prime$ are incompatible $\leftrightarrow$ $\sigma (p)$ and $\sigma (p^\prime )$ are incompatible)
\smallskip\\
(3) $\forall q \in \Q \ \exists p \in \P \ \forall p^\prime \in \P \ (p^\prime \leq p \rightarrow (\sigma (p^\prime )$ and $q$ are compatible in $\Q$)).
\smallskip\\
In (3), we call $p$ a reduction of $q$ to $\P$ with respect to $\sigma$.
\medskip\\
If only (1) and (2) hold, we say that $\sigma$ is an embedding. If $\P \subseteq \Q$ such that the identity is an embedding, then we write $\P \subseteq _\bot \Q$.
\medskip\\
We say that $\P \subseteq \Q$ is completely contained in $\Q$ if $id \upharpoonright \P :\P \rightarrow \Q$ is a complete embedding.
\bigskip\\
Let $\alpha \in Lim$. A finite support (FS) iteration is a sequence $\langle \P _\xi \mid \xi \leq \alpha \rangle$ of partial orders together with a commutative system $\langle \sigma _{\xi \eta} \mid \xi < \eta \leq \alpha \rangle$ of complete embeddings $\sigma _{\xi \eta}:\P _\xi \rightarrow \P _\eta$ such that $\bigcup \{ \sigma _{\xi \eta}[\P _\xi ] \mid \xi < \eta \}=\P _\eta$ for limit $\eta$.
\medskip\\
This is the original definition by Solovay and Tennenbaum in \cite{SolovayTennenbaum}, except that they use Boolean algebras instead of partial orders. Moreover, it is well known that if $\sigma:\P_1 \rightarrow \P_2$ is a complete embedding then there is a $\P_1$-name $\dot \Q$ such that $\P_2$ and $\P_1 \ast \dot \Q$ are forcing equivalent. This leads to the more common definition of FS iterations where conditions are sequences of names. For the exact relationship between the two approaches see Kunen's textbook \cite{Kunen}, chapter VIII \S 5 and exercise K. 
\bigskip\\
An important property of FS iterations is that they preserve the $\kappa$-cc:
\pagebreak\\
{\bf Theorem 2.1}
\smallskip\\
Let $\langle\langle \P _\xi \mid \xi \leq \alpha \rangle ,\langle \sigma _{\xi \eta} \mid \xi < \eta \leq \alpha \rangle \rangle$ be a FS iteration. Assume that all $\P _\xi$ with $\xi < \alpha$ satisfy the $\kappa$-cc. Then $\P _\alpha$ also satisfies the $\kappa$-cc.
\smallskip\\
{\bf Proof:} See the original article by Solovay and Tennenbaum \cite{SolovayTennenbaum} or any standard textbook. $\Box$

\section{Morasses}
A simplified $(\kappa, 1)$-morass is a structure $\frak{M}=\langle \langle \theta _\alpha \mid \alpha \leq \kappa \rangle , \langle \frak{F}_{\alpha\beta}\mid \alpha < \beta \leq \kappa\rangle\rangle$ satisfying the following conditions:
\smallskip\\
(P0) (a) $\theta _0=1$, $\theta _\kappa =\kappa^+$, $\forall \alpha < \kappa \ \ 0<\theta_\alpha < \kappa$.
\smallskip\\
(b) $\frak{F}_{\alpha\beta}$ is a set of order-preserving functions $f:\theta _\alpha \rightarrow \theta _\beta$.
\smallskip\\
(P1) $|\frak{F}_{\alpha\beta}| < \kappa$ for all $\alpha < \beta < \kappa$.
\smallskip\\
(P2) If $\alpha < \beta < \gamma$, then $\frak{F}_{\alpha\gamma}=\{ f \circ g \mid f \in \frak{F}_{\beta\gamma}, g \in \frak{F}_{\alpha\beta}\}$.
\smallskip\\
(P3) If $\alpha < \kappa$, then $\frak{F}_{\alpha, \alpha +1}=\{ id \upharpoonright \theta _\alpha , f_\alpha\}$ where $f_\alpha$ is such that $f_\alpha \upharpoonright \delta = id \upharpoonright \delta$ and $f_\alpha (\delta) \geq \theta _\alpha$ for some $\delta < \theta _\alpha$.
\smallskip\\
(P4) If $\alpha \leq \kappa$ is a limit ordinal, $\beta _1,\beta_2 < \alpha$ and $f_1 \in \frak{F}_{\beta _1\alpha}$, $f_2 \in \frak{F}_{\beta _2\alpha}$, then there are a $\beta _1, \beta _2 < \gamma < \alpha$, $g \in \frak{F}_{\gamma\alpha}$ and  $h_1 \in \frak{F}_{\beta _1\gamma}$, $h_2 \in \frak{F}_{\beta _2\gamma}$ such that $f_1=g \circ h_1$ and $f_2=g \circ h_2$.
\smallskip\\
(P5) For all $\alpha >0$, $\theta _\alpha =\bigcup \{ f[\theta _\beta] \mid \beta < \alpha , f \in \frak{F}_{\beta\alpha}\}$.
\bigskip\\
{\bf Lemma 3.1}
\smallskip\\
Let $\alpha < \beta \leq \kappa$, $\tau _1,\tau _2 < \theta _\alpha$, $f_1,f_2 \in \frak{F}_{\alpha \beta}$ and $f_1(\tau _1)=f_2(\tau _2)$. Then $\tau _1=\tau _2$ and $f_1\upharpoonright \tau _1 = f_2 \upharpoonright \tau _2$.
\smallskip\\
{\bf Proof} by induction over $\beta$: The base case of the induction is $\beta =\alpha +1$. Then the claim follows immediately from (P3). So assume that $\beta = \gamma +1$. Let, by (P2), $f_i=g_i \circ f^\prime_i$ where $f^\prime_i \in \frak{F}_{\alpha\gamma}$, $g_i \in \frak{F}_{\gamma\beta}$. Let $\tau^\prime _i=f^\prime _i(\tau _i)$. It follows like in the base case that $\tau ^\prime_1=\tau^\prime_2$ and $f^\prime_1 \upharpoonright \tau _1=f^\prime _2 \upharpoonright \tau _2$. So, by the induction hypothesis, $\tau _1=\tau _2$ and $f^\prime _1 \upharpoonright \tau _1 =f^\prime _2 \upharpoonright \tau _2$. Hence $f_1\upharpoonright \tau _1 = f_2 \upharpoonright \tau _2$.

Finally, let $\beta \in Lim$. Then there exists by (P4) $\alpha < \gamma < \beta$ and $g \in \frak{F}_{\gamma\beta}$ such that $f_i=g \circ f^\prime_i$, $f^\prime_i \in \frak{F}_{\alpha\gamma}$. So $f^\prime_1(\tau_1)=f^\prime_2(\tau_2)$. Hence $\tau_1=\tau_2$ and  $f^\prime _1 \upharpoonright \tau _1 =f^\prime _2 \upharpoonright \tau _2$ by the induction hypothesis. Therefore $f_1\upharpoonright \tau _1 = f_2 \upharpoonright \tau _2$. $\Box$
\bigskip\\
A simplified morass defines a tree $\langle T , \prec \rangle$:
\medskip\\
Let $T=\{ \langle \alpha , \gamma\rangle \mid \alpha \leq \kappa , \gamma < \theta _\alpha\}$.
\smallskip\\
For $t =\langle \alpha , \nu \rangle \in T$ set $\alpha (t)=\alpha$ and $\nu (t)=\nu$.
\smallskip\\
Let $\langle\alpha , \nu \rangle \prec \langle\beta , \tau\rangle$ iff $\alpha <\beta$ and $f(\nu)=\tau$ for some $f \in \frak{F}_{\alpha \beta}$.
\smallskip\\
If $s:=\langle\alpha , \nu \rangle \prec \langle\beta , \tau\rangle =: t$, $f \in \frak{F}_{\alpha\beta}$ and $f(\nu)=\tau$, then $f \upharpoonright (\nu (s) +1)$ does not depend on $f$ by lemma 3.1. So we may define $\pi _{st}:=f \upharpoonright (\nu (s) +1)$.
\pagebreak\\
{\bf Lemma 3.2}
\smallskip\\
The following hold:
\smallskip\\
(a) $\prec$ is a tree, $ht_T(t)=\alpha (t)$.
\smallskip\\
(b) If $t_0 \prec t_1 \prec t_2$, then $\pi _{t_0t_1}=\pi _{t_1t_2} \circ \pi _{t_0t_1}$.
\smallskip\\
(c) Let $s \prec t$ and $\pi =\pi _{st}$. If $\pi (\nu^\prime)=\tau^\prime$, $s^\prime=\langle \alpha (s), \nu ^\prime\rangle$ and $t^\prime=\langle \alpha (t), \tau^\prime\rangle$, then $s^\prime \prec t^\prime$ and $\pi _{s^\prime t^\prime}=\pi \upharpoonright (\nu^\prime +1)$.
\smallskip\\
(d) Let $\gamma \leq \kappa$, $\gamma \in Lim$. Let $t \in T_\gamma$. Then $\nu (t) +1=\bigcup \{ rng(\pi_{st})\mid s \prec t\}$.
\smallskip\\
{\bf Proof:} (a) First, we prove that $\prec$ is transitive. Let $\langle \alpha , \nu \rangle \prec \langle \beta , \tau \rangle$ be witnessed by $f \in \frak{F}_{\alpha\beta}$ and $\langle \beta , \tau \rangle \prec \langle \gamma , \eta \rangle$ by $g \in \frak{F}_{\alpha\beta}$. Set $h=g \circ f \in \frak{F}_{\alpha\gamma}$ by (P2). Then $h(\nu)=\eta$. So $\langle \alpha , \nu \rangle \prec \langle \gamma, \eta \rangle$. 

Now, let $\langle \alpha , \nu \rangle , \langle \beta , \tau \rangle \prec \langle \gamma , \eta \rangle$ and $\langle \alpha , \nu \rangle \neq \langle \beta , \tau \rangle$. It follows from lemma 3.1 that $\alpha \neq \beta$. Let w.l.o.g. $\alpha < \beta$. Let $\langle \alpha , \nu \rangle \prec \langle \gamma,\eta \rangle$ be witnessed by $f \in \frak{F}_{\alpha\gamma}$. By (P2) choose $g \in \frak{F}_{\beta\gamma}$ and $h \in \frak{F}_{\alpha\beta}$ such that $f=g \circ h$. Then $\langle \alpha , \nu \rangle \prec \langle \beta , h(\nu)\rangle \prec \langle \gamma, \eta \rangle$. However, $h(\nu)=\tau$ by lemma 3.1. Hence $\langle\alpha , \nu \rangle \prec \langle \beta , \tau \rangle$. This proves that $\prec$ is a tree.

Finally, by (P2), for all $t \in T$ there is $s \prec t$ such that $\alpha (s)=\beta$ if $\beta < \alpha (t)$. This shows the second claim.
\smallskip\\
(b) follows immediately from (a) and the definition.
\smallskip\\
(c) Let $s \prec t$ be witnessed by $f \in \frak{F}_{\alpha\beta}$. Then $s^\prime \prec t^\prime$ is also witnessed by $f$ and  $\pi _{s^\prime t^\prime}=\pi \upharpoonright (\nu^\prime +1)$ holds by definition.
\smallskip\\
(d) It suffices to prove $\subseteq$. Let $\nu =\nu(t)$ and $\tau < \nu$. By (P5) choose $\alpha _1,\alpha _2 < \gamma$ and $f_i \in \frak{F}_{\alpha _i \gamma}$ such that $\tau \in rng(f_1)$ and $\nu \in rng(f_2)$. By (P4) choose $\beta$ such that $\alpha _1,\alpha _2 < \beta < \gamma$ and $f^\prime _i \in \frak{F}_{\alpha _i\beta}$, $g \in \frak{F}_{\beta\gamma}$ where $f_i =g \circ f^\prime_i$. Then $\tau,\nu \in rng(g)$. So let $g(\bar\tau)=\tau$ and $g(\bar\nu)=\nu$. Hence $\bar\tau < \bar\nu$, since $g$ is order-preserving. Let $s=\langle \beta , \bar \nu\rangle$. Then $s \prec t$ and $\pi _{st}(\bar \tau)=\tau$. $\Box$
\bigskip\\
{\bf Lemma 3.3}
\smallskip\\
Let $\alpha < \beta \leq \kappa$. Then $id \upharpoonright \theta_\alpha \in \frak{F}_{\alpha\beta}$.
\smallskip\\
{\bf Proof} by induction on $\beta$: The base case of the induction is $\beta=\alpha+1$. Then the claim is part of (P3). So assume that $\beta=\gamma+1$. By the induction hypothesis, $id \upharpoonright \theta_\alpha \in \frak{F}_{\alpha\gamma}$. By (P3), $id \upharpoonright \theta_\gamma \in \frak{F}_{\gamma\beta}$. Hence $id \upharpoonright \theta_\alpha = (id \upharpoonright \theta_\gamma) \circ (id \upharpoonright \theta_\alpha) \in \frak{F}_{\alpha\beta}$ by (P2).

Finally, let $\beta \in Lim$. Assume towards a contradiction that $id \upharpoonright \theta_\alpha \notin \frak{F}_{\alpha\beta}$. Let $f \in \frak{F}_{\alpha\beta}$ be such that $sup(f[\theta_\alpha])$ is minimal. Since $f \neq id \upharpoonright \theta_\alpha$, there are $\nu < \tau \in \theta_\beta$ such that $\nu \notin rng(f)$ and $\tau \in rng(f)$. Let $t =\langle \beta ,\tau\rangle$. By lemma 3.2 (d), there is an $s \prec t$ such that $\nu \in rng(\pi_{st})$. Let $s=\langle \gamma+1,\bar \tau\rangle$ be the minimal such $s$. Let $\nu =\pi_{st}(\bar \nu)$. Furthermore, let $f=f_3 \circ f_2 \circ f_1$ where $f_3 \in \frak{F}_{\gamma+1,\beta}$, $f_2 \in \frak{F}_{\gamma,\gamma+1}$ and $f_1 \in \frak{F}_{\alpha\gamma}$. Then $\pi_{st}=f_3 \upharpoonright \bar \tau +1$. Hence by the minimality of $s$, $f_2 \neq id \upharpoonright \theta_\gamma$, $\bar \nu < \theta_\gamma$ and $\bar \tau \geq \theta_\gamma$. Define $g:= f_3 \circ (id \upharpoonright \theta_\gamma) \circ f_1$. Then $g \in \frak{F}_{\alpha\gamma}$ by (P2) and $rng(g) \subseteq f_3[\theta_\gamma] \subseteq f_3(\bar \tau)=\tau$. Hence $sup(f[\theta_\alpha])$ was not minimal. Contradiction! $\Box$
\pagebreak\\
{\bf Theorem 3.4}
\smallskip\\
(a) If $V=L$, then there is a simplified $(\kappa ,1)$-morass for all regular $\kappa >\omega$.
\smallskip\\
(b) If $\kappa$ is an uncountable regular cardinal such that $\kappa ^+$ is not inaccessible in $L$, then there is a simplified $(\kappa ,1)$-morass.
\smallskip\\
(c) For every regular $\kappa > \omega$, there is a $\kappa$-complete (i.e. every decreasing sequence of length $<$ $\kappa$ has a lower bound) forcing $\P$ satisfying $\kappa ^+$-cc such that $\P \Vdash ($ there is a simplified $(\kappa ,1)$-morass).
\smallskip\\
{\bf Proof:} (a) see Devlin \cite{Devlin}, VIII 2 and 4 or Velleman \cite{Velleman1984}.
\smallskip\\
(b) see Devlin \cite{Devlin}, VIII 4 and exercise 6, or Velleman \cite{Velleman1984b}.
\smallskip\\
(c) see Velleman \cite{Velleman1984}. $\Box$

\section{FS Systems Along Morasses}
Let $\frak{M}$ be a simplified $(\kappa ,1)$-morass. We want to define a generalization of a FS iteration which is not indexed along an ordinal but along $\frak{M}$. One way of doing this is the following definition:
\medskip\\
We say that $\langle\langle \P _\eta \mid \eta \leq \kappa ^+ \rangle ,\langle \sigma _{st} \mid s \prec t \rangle , \langle e_\alpha \mid \alpha < \kappa \rangle\rangle$ is a FS system along $\frak{M}$ if the following conditions hold:
\medskip\\
(FS1) $\langle \P _\eta \mid \eta \leq \kappa ^+  \rangle$ is a sequence of partial orders such that $\P _\eta \subseteq _\bot \P _\nu$ if $\eta \leq \nu$ and $\P _\lambda =\bigcup \{ \P _\eta\mid \eta < \lambda\}$ for $\lambda \in Lim$.
\medskip\\
(FS2) $\langle \sigma _{st} \mid s \prec t \rangle$ is a commutative system of injective embeddings $\sigma _{st}:\P _{\nu (s)+1} \rightarrow \P _{\nu (t)+1}$ such that if $t$ is a limit point in $\prec$, then $\P _{\nu (t) +1} = \bigcup \{ \sigma _{st}[\P _{\nu (s)+1}]\mid s \prec t \}$. 
\medskip\\
(FS3) $e_\alpha : \P _{\theta _{\alpha +1}} \rightarrow \P _{\theta _\alpha}$.
\medskip\\
(FS4) Let $s \prec t$ and $\pi =\pi _{st}$. If $\pi (\nu^\prime)=\tau^\prime$, $s^\prime=\langle \alpha (s), \nu ^\prime\rangle$ and $t^\prime=\langle \alpha (t), \tau^\prime\rangle$, then $\sigma _{st}:\P_{\nu (s)+1} \rightarrow \P _{\nu (t) +1}$ extends $\sigma _{s^\prime t^\prime}:\P _{\nu ^\prime +1}\rightarrow \P _{\tau^\prime +1}$.
\medskip\\
Hence for $f \in \frak{F}_{\alpha\beta}$, we may define $\sigma _f  =\bigcup \{ \sigma _{st} \mid s=\langle \alpha , \nu \rangle, t=\langle \beta , f(\nu)\rangle\}$.
\medskip\\
(FS5) If $\pi _{st}=id \upharpoonright \nu (s)+1$, then $\sigma _{st}=id \upharpoonright \P _{\nu (s) +1}$.
\medskip\\
(FS6)(a) If $\alpha < \kappa$, then $\P _{\theta_\alpha}$ is completely contained in $\P _{\theta _{\alpha +1}}$ in such a way that $e_\alpha(p)$ is a reduction of $p \in \P _{\theta _{\alpha +1}}$. 
\smallskip\\
(b) If $\alpha < \kappa $, then $\sigma _\alpha :=\sigma _{f_\alpha} :\P _{\theta _\alpha} \rightarrow \P _{\theta _{\alpha +1}}$ is a complete embedding such that $e_\alpha(p)$ is a reduction of $p \in \P _{\theta _{\alpha +1}}$.
\medskip\\
(FS7)(a) If $\alpha < \kappa$ and $p \in \P _{\theta _\alpha}$, then $e_\alpha(p)=p$.
\smallskip\\
(b) If $\alpha < \kappa$ and $p \in rng(\sigma _\alpha)$, then $e_\alpha(p)=\sigma ^{-1}_\alpha(p)$.
\bigskip\\
To simplify notation, set $\P := \P _{\kappa ^+}$.
\medskip\\
Unlike in the case of FS iterations, it is unclear how a generic extension with respect to $\P_{\kappa^+}$ can be viewed as being obtained by successive extensions. This would justify to call a FS system along $\frak{M}$ a FS {\it iteration} along $\frak{M}$.
\medskip\\
However, like in the case of FS iterations it is sometimes more convenient to represent $\P$ as a set of functions $p^\ast:\kappa \rightarrow V$ such that $p^\ast(\alpha) \in \P _{\theta _\alpha}$ for all $\alpha < \kappa$.
\bigskip\\
To define such a function $p^\ast$ from $p \in \P$ set recursively
\smallskip

$p_0=p$
\smallskip

$\nu _n(p)=min\{ \eta \mid p_n \in \P _{\eta +1}\}$
\smallskip

$t_n(p)=\langle \kappa ,\nu _n(p)\rangle$
\smallskip

$p^{(n)}(\alpha)=\sigma ^{-1}_{st}(p_n)$ if $s \in T_\alpha$, $s \prec t_n(p)$ and $p_n \in rng(\sigma_{st})$.
\medskip\\
Note that, by lemma 3.2 (a), $s$ is uniquely determined by $\alpha$ and $t_n(p)$. Hence we really define a function. Set

$\gamma _n(p)=min(dom(p^{(n)}))$.
\medskip\\
By (FS2), $\gamma _n(p)$ is a successor ordinal or $0$. Hence, if $\gamma_n(p)\neq 0$, we may define
\smallskip

$p_{n+1}=e_{\gamma_n(p)-1}(p^{(n)}(\gamma _n(p)))$.
\smallskip\\
If $\gamma_n(p)= 0$, we let $p_{n+1}$ be undefined.
\medskip\\ 
Finally, set $p^\ast=\bigcup\{ p^{(n)} \upharpoonright [\gamma _n(p) , \gamma _{n-1}(p)[ \ \mid n \in \omega \}$ where $\gamma _{-1}(p)=\kappa$.
\medskip\\
Note: If $n>0$ and $\alpha \in [\gamma _n(p),\gamma_{n-1}(p)[$, then $p^\ast(\alpha)=\sigma _{s\bar t}^{-1}(p_n)$ where $\bar t=\langle \gamma_n(p)-1,\nu _n(p)\rangle$ because $p^\ast(\alpha)=p^{(n)}(\alpha)=\sigma^{-1}_{st}(p_n)=(\sigma _{\bar tt}\circ \sigma _{s\bar t})^{-1}(p_n)=\sigma_{s\bar t}(p_n)$ where the first two equalities are just the definitions of $p^\ast$ and $p^{(n)}$. For the third equality note that $\bar t \prec t$ since $id \upharpoonright \theta _\alpha \in \frak{F}_{\alpha\beta}$ for all $\alpha < \beta \leq \kappa$ by lemma 3.3. So the equality follows from the commutativity of $\langle \sigma _{st} \mid s \prec t \rangle$. The last equality holds by (FS5).
\medskip\\
It follows from the previous observation that $\langle \gamma _n(p) \mid n \in \omega \rangle$ is decreasing. So the recursive definition above breaks down at some point, i.e. $\gamma _n(p)=0$ for some $n \in \omega$. However, that is good news because of the following.
\medskip\\
The support of $p$ is defined by $supp(p)=\{ \gamma _n(p) \mid n \in \omega\}$. Hence $supp(p)$ is finite.  
\bigskip\\
{\bf Lemma 4.1}
\smallskip\\
If $p^\ast(\alpha)$ and $q^\ast(\alpha)$ are compatible for $\alpha=max(supp(p) \cap supp(q))$, then $p$ and $q$ are compatible.
\smallskip\\
{\bf Proof:} Suppose that $p$ and $q$ are incompatible. Without loss of generality let $\nu := min\{ \eta \mid p \in \P_{\eta +1}\} \leq min \{ \eta \mid q \in \P _{\eta +1} \} =:\tau$. Set $s=\langle \kappa,\nu\rangle$ and $t=\langle \kappa, \tau \rangle$. Let $t^\prime \prec t$ be minimal such that $\nu \in rng(\pi _{t^\prime t})$ and $p,q \in rng(\sigma _{t^\prime t})$. By (FS2), $t^\prime \in T_{\alpha _0 +1}$ for some $\alpha < \kappa$. Let $\pi _{t^\prime t}(\nu ^\prime)=\nu$ and $s^\prime =\langle \alpha +1, \nu ^\prime\rangle$. Let $\bar s$, $\bar t$ be the direct predecessors of $s^\prime$ and $t^\prime$ in $\prec$. Set $p^\prime = \sigma _{s^\prime s}^{-1}(p)$, $q^\prime=\sigma ^{-1}_{t^\prime t}(q)$. Then $p^\prime =p^\ast(\alpha _0+1)$, $q^\prime =q^\ast(\alpha _0+1)$ by the definition of $p^\ast$. Moreover, $p^\prime$ and $q^\prime$ are not compatible, because if $r \leq p^\prime , q^\prime$, then $\sigma _{t^\prime t}(r) \leq p,q$ by (FS2). Now, we consider several cases.
\medskip\\
{\it Case 1}: $\nu ^\prime \notin rng(\pi _{\bar t t^\prime})$
\smallskip\\
Then $\pi _{\bar s s^\prime}=id \upharpoonright \nu (\bar s)+1$ and $\sigma _{\bar s s^\prime}=id \upharpoonright \P _{\nu (\bar s)+1}$ by the minimality of $\alpha _0$. Moreover, $\bar p:=p^\prime$ and $\bar q:=e_\alpha (q^\prime)$ are not compatible, because if $r \leq p^\prime, e_\alpha (q^\prime)$, then there is $u \leq r,q^\prime,p^\prime$ by (FS6)(a). There is no difference between compatibility in $\P _{\theta _{\alpha +1}}$ and in $\P _{\nu (t^\prime)+1}$ by (FS1). Finally, note that $\bar p=p^\ast(\alpha _0)$ and $\bar q =q^\ast(\alpha _0)$ by the definition of $p^\ast$ and (FS7).
\medskip\\
{\it Case 2}:  $\nu ^\prime \in rng(\pi _{\bar t t^\prime})$ and $\pi _{\bar ss^\prime}=id \upharpoonright \nu (\bar s)+1$
\smallskip\\
Then $\pi _{\bar t t^\prime} \neq id \upharpoonright \nu (\bar t)+1$ by the minimality of $\alpha _0$ and $\bar p:=p^\prime$ and $\bar q:=e_\alpha(q^\prime)$ are not compatible (like in case 1). However, $\bar p=p^\ast(\alpha _0)$ and $\bar q =q^\ast(\alpha _0)$ by the definition of $p^\ast$ and (FS7).
\medskip\\
{\it Case 3}: $\nu ^\prime \in rng(\pi _{\bar t t^\prime})$, $\pi _{\bar ss^\prime}\neq id \upharpoonright \nu (\bar s)+1$ and $\alpha_0 +1 \notin supp(p)$
\smallskip\\
Then $\pi _{\bar t t^\prime} \neq id \upharpoonright \nu (\bar t)+1$ by the minimality of $\alpha _0$. Set $\bar p:=\sigma ^{-1}_{\bar s s^\prime}(p^\prime)$ and $\bar q = e_\alpha (q^\prime)$. Then $\bar p$ and $\bar q$ are not compatible, because if $r \leq \bar p, \bar q$, then there is $u \leq \sigma _\alpha (r),q^\prime , p^\prime$ by (FS6)(b). However, $\bar p=p^\ast(\alpha _0)$ and $\bar q =q^\ast(\alpha _0)$ by the definition of $p^\ast$ and (FS7).
\medskip\\
{\it Case 4}: $\nu ^\prime \in rng(\pi _{\bar t t^\prime})$, $\pi _{\bar ss^\prime}\neq id \upharpoonright \nu (\bar s)+1$ and $\alpha_0 +1 \notin supp(q)$ 
\smallskip\\
Then $\pi _{\bar t t^\prime} \neq id \upharpoonright \nu (\bar t)+1$. Set $\bar q:=\sigma ^{-1}_{\bar s s^\prime}(q^\prime)$ and $\bar p = e_\alpha (p^\prime)$. Then $\bar q$ and $\bar p$ are not compatible, because if $r \leq \bar p, \bar q$, then there is $u \leq \sigma _\alpha (r),p^\prime , q^\prime$ by (FS6)(b).
\medskip\\
{\it Case 5}: $\alpha_0 +1 \in supp(p) \cap supp(q)$
\smallskip\\
Then $\alpha_0 +1 = max(supp(p) \cap supp(q))$, since $\alpha_0+1 \geq max(supp(q))$ because by definition $q \in rng(\sigma _{rt})$ where $r \prec t$ and $r \in T_{max(supp(q))}$. However, $p^\prime=p^\ast(\alpha_0+1)$, $q^\prime=q^\ast(\alpha_0+1)$ are not compatible. Contradiction.
\medskip\\
So in case 5 we are finished. If we are in cases 1 - 4, we define recursively $\alpha _{n+1}$ from $p^\ast (\alpha_n)$ and $q^\ast(\alpha _n)$ in the same way as we defined $\alpha_0$ from $p$ and $q$. Like in the previous proof that $\langle \gamma _n(p)\mid n \in \omega \rangle$ is decreasing, we see that $\langle \alpha _n \mid n \in \omega \rangle$ is decreasing. Hence the recursion breaks off, we end up in case 5 and get the desired contradiction. $\Box$
\bigskip\\
{\bf Theorem 4.2}
\smallskip\\
Let $\mu ,\kappa >\omega$ be cardinals, $\kappa$ regular. Let  $\langle\langle \P _\eta \mid \eta \leq \kappa ^+ \rangle ,\langle \sigma _{st} \mid s \prec t \rangle ,\langle e_\alpha \mid \alpha < \kappa \rangle \rangle$ be a FS system along a $(\kappa ,1)$-morass $\frak{M}$. Assume that all $\P _\eta$ with $\eta < \kappa$ satisfy the $\mu$-cc. Then $\P _{\kappa ^+}$ also does.
\smallskip\\
{\bf Proof:} Let $A \subseteq \P_{\kappa^+}$ be a set of size $\mu$. Assume by the $\Delta$-system lemma that $\{ supp(p) \mid p \in A\}$ forms a $\Delta$-system with root $\Delta$. Set $\alpha = max(\Delta)$. Then $\P _{\theta _\alpha}$ satisfies the $\mu$-cc by the hypothesis of the lemma. So there are $p \neq q \in A$ such that $p^\ast(\alpha)$ and $q^\ast(\alpha)$ are compatible. Hence $p$ and $q$ are compatible by the previous lemma. $\Box$

\section{A CCC Forcing That Adds An $\omega_2$-Suslin Tree}
As an application, we construct along an $(\omega _1,1)$-morass a ccc forcing $\P$ that adds an $\omega_2$-Suslin tree.
\smallskip\\
The natural forcing to do this with finite conditions is Tennenbaum's forcing (see \cite{Tennenbaum}): Define $P(\theta)$ as the set of all finite trees $p=\langle x_p,<_p\rangle$, $x_p \subseteq \theta$, such that $\alpha < \beta$ if $\alpha <_p\beta$.
\smallskip\\
Set $p \leq q$ iff $x_p \supset x_q$ and $<_q=<_p \cap x^2_q$.
\smallskip\\
For $\theta =\omega_1$, $P(\theta)$ is Tennenbaum's forcing to add an $\omega_1$-Suslin tree which satisfies ccc.
\smallskip\\
However, if $\theta > \omega_1+1$, then
\smallskip

$A=\{ p \in P(\theta) \mid x_p=\{ \alpha , \alpha +1,\alpha +2, \omega _1, \omega_1+1\} , \alpha < \omega _1,$
$$\alpha <_p \alpha +1<_p \omega _1,\alpha <_p \alpha +2<_p\omega_1+1 , \alpha +1 \not <_p\alpha +2\}$$
is an antichain of size $\omega_1$.
\smallskip\\
So $P(\theta)$ does not satisfy the ccc and in order to thin it out so that it obtains ccc, we have to restrict the possible values of the infima in our trees. 
\smallskip\\
Let $\pi :\bar \theta \rightarrow \theta$ be a order-preserving map. Then $\pi :\bar\theta \rightarrow \theta$ induces maps $\pi : \bar\theta ^2 \rightarrow \theta ^2$ and $\pi :P(\bar \theta) \rightarrow P(\theta)$ in the obvious way:
$$\pi : \bar\theta ^2 \rightarrow \theta ^2,\quad \langle \alpha,\beta\rangle \mapsto \langle \pi (\alpha ),\pi (\beta) \rangle$$
$$\pi :P(\bar \theta) \rightarrow P(\theta),\quad \langle x_p,<_p\rangle \mapsto \langle \pi [x_p],\pi[<_p]\rangle .$$
If $p \in P(\theta)$, then set
$$\pi ^{-1}[p]:=\langle \pi ^{-1}[x_p \cap rng(\pi)], \pi ^{-1}[<_p \cap rng(\pi)]\rangle.$$
It is easily seen that then $\pi ^{-1}[p] \in P(\bar\theta)$.
\medskip\\
Now, let us assume that we restrict the allowed values of the infimum $i_p(\alpha, \beta)$ of $\alpha , \beta \in x_p$ in the tree $p \in P(\theta)$ to a set $F(\alpha , \beta)$. For $\delta < \theta$, we want to find a reduction of $p \in P(\theta)$ with respect to $id \upharpoonright \delta$. Let us look for example at a $p$ with $x_p=\{ \alpha , \beta\}$ and $<_p=\{ \langle \alpha , \beta \rangle \}$ such that $\alpha < \delta < \beta$. Then we cannot just take $(id \upharpoonright \delta)^{-1}$ as reduction because we could extend it to a condition $q$ such that $i_q(\alpha , \gamma)$ exists for some $\gamma \in \delta$. However, $i_q(\alpha , \gamma)$ could have any value in $F(\alpha , \gamma)$, while in a common extension $r$ of $p$ and $q$ we have $i_r(\gamma , \beta)=i_q(\alpha , \gamma)$ and $i_r(\alpha , \beta)$ has to be an element of $F(\alpha , \beta)$. We can solve this problem by taking $s$ with $x_s=\{ \alpha , \beta^\prime\}$ and $<_s=\{ \langle \alpha , \beta ^\prime\rangle\}$ as reduction for some $\beta^\prime$ with $F(\alpha , \beta^\prime)=F(\alpha , \beta)$. This leads to the following definition in which the $F(\alpha , \beta)$ are not needed anymore. But they could be introduced as the ranges of the morass maps. A similar problem arises in Baumgartner's and Shelah's forcing to add a a thin-very tall superatomic Boolean algebra \cite{BaumgartnerShelah}. They explicitly define a function $F$ like above, which they obtain by historic forcing.    
\medskip\\
We define our FS system by induction on $\beta \leq \omega_1$.
\medskip\\
{\it Base Case}: $\beta =0$
\medskip\\
Then we need to define only $\P _1$. Set $\P _1:=P(1)$.
\pagebreak\\
{\it Successor Case}: $\beta = \alpha +1$
\medskip\\
We first define $\P _{\theta _\beta}$. To do so, let
\medskip

$\P ^\prime _{\theta _\beta} :=\{ \langle x_p \cup x_{f_\alpha(p)}, <_p \cup <_{f_\alpha(p)}\rangle \mid p \in \P_{\theta_\alpha}\}$
$$\cup \{ \langle x_p \cup x_{f_\alpha(p)}, tc(<_p \cup <_{f_\alpha(p)} \cup\{ \langle  \eta , min\{ \gamma \in [\theta _\alpha , \theta_{\alpha +1}[ \mid \gamma \leq _{f_\alpha (p)} f_\alpha (\eta)\} \rangle\}\rangle $$ $$\mid p \in \P_{\theta_\alpha} , \eta \in x_p, \eta < f_\alpha (\eta)\} .$$
In this definition, $tc(x)$ denotes the transitive closure of the binary relation $x$. As we know from Tennenbaum's original proof, every element of $\P^\prime_{\theta_\beta}$ is an element of $P(\theta_\beta)$ which extends $p$ and $f_\alpha(p)$. This is easily seen.
\medskip\\  
Now, define
$$\P _{\theta _\beta} := \{ p \in P(\theta _\beta) \mid r \leq p \hbox{ for some }r \in \P ^\prime _{\theta_\beta}\} .$$
For $t \in T_\beta$ set $\P _{\nu (t)+1}=P(\nu (t)+1) \cap \P _{\theta _\beta}$ and $\P _\lambda =\bigcup \{ \P _\eta \mid \eta < \lambda\}$ for $\lambda \in Lim$. Let $\sigma _{st}:\P _{\nu (s)+1} \rightarrow \P_{\nu(t)+1}, p \mapsto \pi_{st}(p)$. 
\smallskip\\
We still need to define $e_\alpha$. If $p \in rng(\sigma _\alpha)$, then set $e_\alpha (p)=\sigma^{-1}_\alpha(p)$. If $p \in \P _{\theta _\alpha}$, then set $e_\alpha(p)=p$. Finally, if $p \notin rng(\sigma_\alpha) \cup \P _{\theta_\alpha}$, then pick an $r \in \P^\prime_{\theta_\beta}$ such that $r \leq p$ and set $e_\alpha(p)=f^{-1}_\alpha[r]$.
\medskip\\
{\it Limit Case}: $\beta \in Lim$
\medskip\\
Then everything is already uniquely determined by (FS1) and (FS2). That is, for $t \in T_\beta$ set $\P _{\nu(t)+1}=\bigcup \{ \sigma _{st}[\P_{\nu (s)+1}] \mid s \prec t \}$ and $\P _\lambda =\bigcup \{ \P _\eta \mid \eta < \lambda\}$ for $\lambda \in Lim$. Let $\sigma _{st}:\P _{\nu (s)+1} \rightarrow \P_{\nu(t)+1}, p \mapsto \pi_{st}(p)$.
\bigskip\\
{\bf Lemma 5.1}
\smallskip\\
$\P$ satisfies the ccc.
\smallskip\\
{\bf Proof:} Since all $P(\theta)$ for $\theta < \omega_1$ have size $\leq \omega$, it suffices by theorem 4.2 to show that $\langle\langle \P _\eta \mid \eta \leq \kappa ^+ \rangle ,\langle \sigma _{st} \mid s \prec t \rangle ,\langle e_\alpha \mid \alpha < \kappa \rangle \rangle$ is a FS system along the morass.
\smallskip\\
Most conditions of the definition of a FS system are clear. We only prove (FS6). Let $p \in \P _{\theta _\beta}$ and $\beta = \alpha +1$. We may assume that $p \in \P ^\prime _{\theta _\beta}$, because by definition $\P _{\theta _\beta}$ is dense in $\P ^\prime _{\theta _\beta}$.  We have to prove that $\sigma_\alpha ^{-1}[p]$ is a reduction of $p$ with respect to $\sigma _\alpha$ and $id \upharpoonright \P _{\theta _\alpha}$. To do so for $\sigma _\alpha$, let $q \leq \sigma_\alpha ^{-1}[p]=:s$. We have to find an $r \leq p,\sigma _\alpha(q)$ such that $r \in \P _{\theta _\beta}$. We consider two cases. If $p$ is of the form $ \langle x_s \cup x_{f_\alpha(s)}, <_s \cup <_{f_\alpha(s)}\rangle$, then define $r:= \langle x_q \cup x_{f_\alpha(q)}, <_p \cup <_{f_\alpha(q)}\rangle$. It is easily seen that this is an extension of $p$ and $\sigma _\alpha(q)$. If $p$ is of the form $$\langle x_s \cup x_{f_\alpha(s)}, tc(<_s \cup <_{f_\alpha(s)} \cup\{ \langle  \eta , min\{ \gamma \in [\theta _\alpha , \theta_{\alpha +1}[ \mid \gamma \leq _{f_\alpha (s)} f_\alpha (\eta)\} \rangle\}\rangle$$ for some $\eta \in x_s$, then define $r$ as
$$\langle x_q \cup x_{f_\alpha(q)}, tc(<_q \cup <_{f_\alpha(q)} \cup\{ \langle  \eta , min\{ \gamma \in [\theta _\alpha , \theta_{\alpha +1}[ \mid \gamma \leq _{f_\alpha (q)} f_\alpha (\eta)\} \rangle\}\rangle.$$  
Again, it is easily seen that this is an extension of $p$ and $\sigma _\alpha(q)$. That proves that $\sigma_\alpha ^{-1}[p]$ is a reduction of $p$ with respect to $\sigma _\alpha$. The proof that $\sigma_\alpha ^{-1}[p]$ is a reduction of $p$ with respect to $id \upharpoonright \P _{\theta _\alpha}$ is completely analogous. $\Box$ 
\bigskip\\
{\bf Lemma 5.2}
\smallskip\\
If $\gamma_0(p)=\gamma_0(q)$, $p^\ast(\gamma_0(p))=q^\ast(\gamma_0(q))$, $\pi:p \cong q$ and $\alpha \leq \pi (\alpha)$, then there exists an $r \leq p,q$ such that $\langle \alpha , \pi (\alpha )\rangle \in \leq_r$.
\smallskip\\
{\bf Proof:} Let $p$ and $q$ be as in the hypothesis of the lemma. We prove by induction over $\eta \in [\gamma_0(p), \omega_1]$ that if $\pi:p^\ast(\eta) \cong q^\ast (\eta)$ (where $p^\ast (\omega_1):=p$) and $\alpha \leq \pi (\alpha)$, then there exists an $r \leq p^\ast(\eta),q^\ast(\eta)$ such that $\langle \alpha , \pi (\alpha)\rangle \in \leq _r$.
\smallskip\\
{\it Base Case}: $\eta = \gamma_0(p)=\gamma_0(q)$
\smallskip\\
In this case the claim is trivial because $p^\ast (\eta)=q^\ast(\eta)$.
\smallskip\\
{\it Successor Case}: $\eta =\gamma +1$
\smallskip\\
Let $\pi:p^\ast (\eta) \cong q^\ast (\eta )$ and $\alpha \leq \pi (\alpha)$. Let $\sigma_p:p^\ast(\gamma) \cong p^\ast(\eta)$, $\sigma_q:q^\ast(\gamma)\cong q^\ast(\eta)$ and $\sigma _p(\bar \alpha_p)=\alpha$, $\sigma _q(\bar \alpha _q)=\pi(\alpha)$. By the induction hypothesis, there is an $s \leq p^\ast(\gamma),q^\ast(\gamma)$ such that $\langle \bar \alpha _q,\bar \alpha _p\rangle \in \leq _s$ or  $\langle \bar \alpha _p,\bar \alpha _q\rangle \in \leq _s$. Let $\bar \alpha := max\{ \bar \alpha _p,\bar \alpha _q\}$. Now, we consider two cases. If $\bar \alpha < f_\gamma (\bar\alpha)$, we define $r$ as
$$\langle x_s \cup x_{f_\alpha(s)}, tc(<_s \cup <_{f_\alpha(s)} \cup\{ \langle  \bar\alpha , min\{ \beta \in [\theta _\gamma , \theta_{\gamma +1}[ \mid \beta \leq _{f_\alpha (s)} f_\gamma (\bar\alpha)\} \rangle\}\rangle.$$
If $\bar\alpha = f_\gamma(\bar\alpha)$, then we define
$$r:=  \langle x_s \cup x_{f_\alpha(s)}, <_s \cup <_{f_\alpha(s)}\rangle.$$ 
In both cases, it is easily seen that $r \leq p^\ast(\eta),q^\ast(\eta)$ and $\langle \alpha , \pi (\alpha )\rangle \in \leq_r$.
\smallskip\\
{\it Limit Case}: $\eta \in Lim$
\smallskip\\
By (FS1) and (FS2), there are a $t \in T_\eta$ and an $s \prec t$ such that $p^\ast(\eta),q^\ast (\eta) \in rng(\sigma_{st})$. Let $s \in T_\gamma$, $\sigma _{st}(\bar\alpha)=\alpha$ and $\sigma _{st} \circ \bar \pi=\pi \circ \sigma _{st}$. Then $\sigma _{st}(p^\ast(\eta))=p^\ast (\gamma)$ and  $\sigma _{st}(q^\ast(\eta))=q^\ast (\gamma)$. Moreover, by the induction hypothesis, there is a $\bar r\leq p^\ast (\gamma),q^\ast(\gamma)$ such that $\langle \bar \alpha , \bar \pi (\bar \alpha)\rangle \in \leq _{\bar r}$. Set $r:=\sigma_{st}(\bar r)$. Then $r$ is as desired. $\Box$ 
\bigskip\\
{\bf Lemma 5.3}
\smallskip\\
If $\gamma_0(p)=\gamma_0(q)$, $p^\ast(\gamma_0(p))=q^\ast(\gamma_0(q))$, $\pi:p \cong q$ and $\alpha \leq \pi (\alpha)$, then there exists an $r \leq p,q$ such that $\langle \alpha , \pi (\alpha )\rangle \notin <_r$.
\smallskip\\
{\bf Proof:} Basically the proof proceeds like the proof of lemma 5.3. However, in the successor case, we always use common extensions of the form $ \langle x_p \cup x_{f_\gamma(p)}, <_p \cup <_{f_\gamma(p)}\rangle$. $\Box$
\bigskip\\
{\bf Theorem 5.4}
\smallskip\\
If there is a simplified $(\omega_1,1)$-morass, then there is a ccc forcing that adds an $\omega_2$-Suslin tree.
\smallskip\\
{\bf Proof:} We show that $\P$ forces an $\omega_2$-Suslin tree. To do so, we prove that the generic tree has neither an antichain nor a chain of size $\omega_2$. First, assume towards a contradiction that there is an antichain of size $\omega_2$. Then there is a $p \in \P$ and by ccc of $\P$ a sequence $\langle \dot x_i \mid i \in \omega_2\rangle$ such that $p \Vdash (\{ \dot x_i \mid i \in \omega_2 \}$ is an antichain). Let $\langle \alpha _i \mid i \in \omega_2\rangle$ and $\langle p_i \mid i \in \omega_2\rangle$ be such that $p_i \leq p$ for all $i \in \omega_2$ and $p_i \Vdash (\dot x_i=\check \alpha_i \wedge \dot x_i \in \check x_{p_i})$. Since $card(\P_{\omega_1})=\omega_1$, there is $q \in \P_{\omega_1}$, $\eta \in \omega_1$ and a subset $X\subseteq \omega_2$ of size $\omega_2$ such that $\gamma_0(p_i)=\eta$ and $p_i^\ast(\gamma_0(p_i))=q$ for all $i \in X$. Hence all $p_i$ with $i \in X$ are isomorphic. Since $x_q$ is finite, there are $i \neq j \in X$ such that $\pi(\alpha_i)=\alpha_j$ and $\alpha_i \leq \alpha_j$ where $\pi :p_i \cong p_j$. By lemma 5.2, there exists an $r \leq p_i,p_j$ such that $\langle \alpha _i,\alpha _j\rangle \in \leq_r$. Hence $r \Vdash (\alpha_i$ and $\alpha_j$ are comparable). That contradicts the definition of $p$. The proof that there is no chain of size $\omega_2$ works the same using lemma 5.3 instead of lemma 5.2.   $\Box$

\bibliographystyle{plain}
\bibliography{biblio}
{\sc Mathematisches Institut, Universit\"at Bonn, Beringstrasse 1, 53115 Bonn, Germany}
\end{document}